\documentclass[11pt]{amsart}
\usepackage{amssymb,stmaryrd}
\usepackage{mathrsfs}
\usepackage{enumitem}  
\usepackage{xcolor}
\usepackage{booktabs}
\usepackage{url}
\usepackage{bm}
\newtheorem{theorem}{Theorem}[section]

\newtheorem{lemma}[theorem]{Lemma}

\newtheorem{proposition}[theorem]{Proposition}
\newtheorem{maintheorem}{Theorem}

{\theoremstyle{definition}

}

\DeclareMathOperator {\M}{M}

\begin{document}
\title{Primitive almost simple IBIS groups with sporadic socle}
\author{Melissa Lee}
\address{School of Mathematics, Monash University, Clayton 3800, Australia} 
\email{melissa.lee@monash.edu}

\date{\today}

\maketitle

\begin{abstract}
An irredundant base $B$ for a permutation group $G\leq \mathrm{Sym}(\Omega)$ is an ordered subset of $\Omega$ with trivial stabiliser such that no base point is fixed by the stabiliser of its predecessors.  Groups whose irredundant bases all have the same size are termed Irredundant Bases of Invariant Size (IBIS) groups, and were introduced by Cameron and Fon-Der-Flaass.  In this paper, we contribute to the classification of primitive IBIS groups by classifying those that are almost simple with sporadic socle.
\end{abstract}

Let $G \leq \mathrm{Sym}(\Omega)$ be a permutation group.  A \textit{base} for $G$ is a subset $B$ of $\Omega$ with trivial pointwise stabiliser.  The cardinality of the smallest base for $G$ is called the \textit{base size} and denoted $b(G)$. Since their introduction by Sims in the 1970s, bases have been of enduring research interest because of their applications in computational group theory. 

An ordered base $B=(b_1, b_2, \dots, b_t)$ for $G$ is \textit{irredundant} if, for each $i$ with $1\leq i\leq t-1$, there is an element of $G$ that fixes $b_1, \dots b_i$, but does not fix $b_{i+1}$.
We say that $G$ is an Irredundant Bases of Invariant Size (IBIS) group if all of its irredundant bases are of the same size.  Note that all such irredundant bases must have size equal to the base size of $G$.
Amongst the known examples of IBIS groups are $A_n$ and $S_n$ acting naturally on $n$ points, as well as the general linear group $G=\mathrm{GL}(V)$ acting on a finite dimensional vector space $V$, where any spanning set forms a base for $G$, but a spanning set is an irredundant base if and only if it is a basis.

A reordered irredundant base is not irredundant in general.  Indeed, Cameron and Fon-Der-Flaass \cite{cameron1995bases}, proved that $G$ is an IBIS group if and only if any reordering of an irredundant base for $G$ is irredundant, and the set of irredundant bases for $G$ forms the set of bases of a matroid.  

 Lucchini et al. \cite{LMM} recently made a major breakthrough in the study of IBIS groups by proving that a primitive IBIS group can belong to one of only three of the O'Nan-Scott types, namely almost simple, affine or diagonal type. They further showed that the primitive IBIS groups of diagonal type all belong to the infinite family $\{\mathrm{PSL}_2(2^f)\times \mathrm{PSL}_2(2^f) \mid f\geq 2\}$, where each group has degree $|\mathrm{PSL}_2(2^f)| = 2^f(2^{2f}-1)$.
 Lee and Spiga classified the primitive almost simple IBIS groups with alternating socle \cite{LS}. Namely, they show that apart from the aforementioned natural actions of $S_n$ and $A_n$, there are 14 further examples and in each case, the socle $A_n$ has $n\leq 8$.

In the present paper, we continue the classification of primitive almost simple IBIS groups by classifying those with a sporadic socle. 

\begin{maintheorem}
\label{mainthm}
An almost simple primitive group $G$ with sporadic socle is IBIS if and only if $(G,H,b(G))$ is one of $(\M_{11},A_6.2,4)$, $(\M_{12},\M_{11},5)$, $(\M_{22},\mathrm{L}_3(4),5)$, $(\M_{23},\M_{22},6)$ or $(\M_{24},\M_{23},7)$.
\end{maintheorem}

The methods used to prove Theorem \ref{mainthm} rely heavily on computation in {\sf{GAP}} and {\sc{Magma}}  \cite{Magma,GAP4}. Code to reproduce the results in this paper are available on the author's website \cite{website}.
\section{Proof of Theorem \ref{mainthm}}
\label{sec:proof}

The following two results show that we need only consider groups of base size at least three in order to prove Theorem \ref{mainthm}, and give an explicit list of cases that we need to investigate.

\begin{lemma}[{\cite[Lemma 2.3]{LMM}}]
\label{no_b(G)=2}
There are no IBIS groups $G$ with non-abelian socle and $b(G)=2$. 
\end{lemma}

\begin{proposition}[{\cite[Theorem 1]{MR2684423}}]
\label{b(G)>2}
The almost simple primitive groups with sporadic socle and $b(G) \geq 3$ are known.
\end{proposition}

Let $G$ be an almost simple group with sporadic socle,  acting primitively on the cosets of a maximal subgroup $H$ of $G$. In general, we expect IBIS groups to be rare, so the majority of the techniques that we employ to prove Theorem \ref{mainthm} involve showing that there exists an irredundant base of $G$ of size larger than $b(G)$, so $G$ is not IBIS. 
This may be achieved in a variety of ways, some of which we now set out.  The labellings of these techniques will allow us to be more concise later on.

\begin{itemize}
\item[(T1)] We consider $G$ as a permutation group on cosets of $H$, and we construct a partial irredundant base of size $b(G)$ by random search.
\item[(T2)] We observe that the stabiliser in $G$ of a coset $Hx$ is $H^x$. Therefore,  we find a set of conjugates $\{H^{x_1}, H^{x_2},\dots H^{x_{b(G)}}\}$ of $H$ by random search such that, writing $H_i = \cap_{j=1}^i H^{x_i}$, we have a descending chain of subgroups $G>H_1> H_2> \dots > H_{b(G)} \neq 1$. This implies that $\{Hx_1, Hx_2, \dots Hx_{b(G)}\}$ is a partial irredundant base for $G$.
\item[(T3)] We choose a subgroup $K<G$, and show that $K$ has an irredundant base of size larger than $b(G)$ in its action on the cosets of $H$.
\end{itemize}

 We now embark on the proof of Theorem \ref{mainthm}.
\subsection{$\bm{G_0\in  \{\text{M}_{11},\text{M}_{12}, \text{M}_{22}, \text{M}_{23},\text{M}_{24},\text{J}_1, \text{J}_2,\text{J}_3,\text{HS},\text{McL}\}}$}

In all of these cases,  we are able to construct the permutation representation of $G$ acting on cosets of $H$ in {\sc Magma} by first constructing $G$ as a matrix group, finding its maximal subgroups and using the {\tt LMGCosetImage} function.  We prove most groups are not IBIS by applying (T1).
The remaining cases are when $(G,H,b(G))$ is one of $(\M_{11},A_6.2,4)$, $(\M_{12},\M_{11},5)$, $(\M_{22},\text{L}_3(4),5)$, $(\M_{23},\M_{22},6)$ or $(\M_{24},\M_{23},7)$; we prove they are IBIS by enumerating a set of orbit representatives on $b(G)$-tuples, and showing that each is either an irredundant base, or is not a base.

\subsection{$\bm{G_0\in  \{ \text{Co}_3, \text{Co}_2,\text{He}, \text{Suz}, \text{Fi}_{22}, \text{Fi}_{23}, \text{Fi}_{24}', \text{Ru}, \text{O'N}, \text{Co}_1,
 \text{HN}\}$}}
Here at least one of the permutation representations equivalent to $G$ acting on cosets of $H$ is not readily available.  In most cases however, we are able to construct $H$ using the {\tt AtlasSubgroup} command as part of the AtlasRep package of {\sf GAP} \cite{atlasrep,GAP4}.  
We then apply (T2) in each case and deduce that $G$ is not IBIS.
\subsection{$\bm{G_0=\text{Ly}}$}
\label{Ly}
In view of Proposition \ref{b(G)>2}, we have $H\in \{G_2(5), 3.\text{McL}:2\}$, where $b(G)=3$ in each case.

Let $H=G_2(5)$. We appeal to (T3) by considering the centraliser $K=C_G(g_5)$ of an element $g_5\in H$ of order 5.  We find $H$ has two classes of elements of order 5 and,  by the centraliser orders, we deduce that both classes are contained in the class labelled 5A in $G$. Now $K=C_G(g_5)$ is maximal in $G$ \cite[p. 174]{ATLAS} and, replacing $g_5$ with a conjugate if necessary,  a construction of $K$ is available as a straight line program in the Online Atlas \cite{onlineATLAS}.
Without loss, we choose $g_5$ to be in class 5A in $H$, so that $|C_H(g_5)| = 2^3.3.5^6$, and find candidates for $C_H(g_5)<K$ using the \texttt{Subgroups} command in {\sc{Magma}}.  In order to construct an partial irredundant base for $K$ of length 3, we first fix the trivial coset $H$ so that we now must find a partial irredundant base for $K_2=C_H(g_5)$ of size 2. We then randomly search $K$ for two elements $x_1, x_2$ such that the following conditions hold: 
\begin{enumerate}
\item $x_1x_2^{-1} \notin H$, so that $Hx_1$ and $Hx_2$ are distinct.
\item $K_2> K_2 \cap K_2^{x_1} >1$, and 
\item  $K_2 \cap K_2^{x_1}> K_2 \cap K_2^{x_1} \cap K_2^{x_2} >1$ so that $\{H,Hx_1,Hx_2\}$ forms a partial irredundant base for $K$ by (T2).
\end{enumerate}
Choosing $K$, $K_2$ in this way allows us to verify each of these conditions in $K$ rather than $G$, which significantly reduces the computational resources required. We find an appropriate $x_1, x_2$ via random search and therefore prove that $G$ is not IBIS.

If instead $H = 3.\text{McL}:2$, we repeat the same process, instead considering the centraliser of an involution in class 2A of $H$. We deduce that $G$ is not IBIS.
\subsection{$\bm{G=\text{J}_4}$}
By Proposition \ref{b(G)>2}, we only need consider the action of $G$ on right cosets of one of its maximal subgroups $H \in \{2^{11}:\M_{24}, 2^{1+12}:3.\M_{22}.2, 2^{10}:\mathrm{L}_5(2)\}$, where $b(G)=3$ in each case.

Suppose $H=2^{11}:\M_{24}$.  Now $H$ is the stabiliser of a vector in the irreducible 112-dimensional representation of $G$ over $\mathbb{F}_2$ \cite[p.190]{ATLAS}. Generators for $G$ and $H$ in this representation, along with a vector $v$ fixed by $H$, are given in the Online Atlas \cite{onlineATLAS}.  The action of $G$ on cosets of $H$ is 
then permutationally isomorphic to the action of $G$ on $v^G$.  We consider a subgroup $K \cong \M_{24}$ of $H$ and find that $v$ along with two further vectors in $v^G$ obtained by random search
have stabiliser order chain $[|M_{24}|,16,2]$ in $K$. The result follows applying (T3).

Now suppose $H=2^{1+12}:3.\M_{22}.2$. We employ technique (T3) with the centraliser of an element $g_5\in H$ of order 5.  There is a single class of elements of order 5 in $G$ and so we have $|C_G(g_5)| = 2^6.3.5.7$, while $|C_H(g_5)| = 2^6.3.5$.  Two elements $x_1,x_2 \in C_G(g_5)$ are representatives for the same coset in $G/H$ if and only if $x_1x_2^{-1} \in C_H(g_5)$. We find appropriate $x_1, x_2$ via random search such that by (T2) we can deduce that $H, Hx_1$ and $Hx_2$ form a partial irredundant base for $C_G(g_5)$, so $G$ is not IBIS.

Finally, suppose $H=2^{10}:\mathrm{L}_5(2)$.  We apply (T3) with a subgroup $K \cong L_5(2)$ of $H$. The construction of such a subgroup $K$ is available in the Online Atlas, and we find that $K$ fixes a non-zero vector $w \in V$.  As before, we consider the action of $G$ on $w^G$.  By random search, we find $h_1, h_2 \in G$ such that $H=G_{w}>G_{w,w^{h_1}}>G_{w,w^{h_1}, w^{h_2}}>1$, so $G$ is not IBIS.

%
%
%
\subsection{$\bm{G=\text{Th}}$}
By Lemma \ref{no_b(G)=2} and Proposition \ref{b(G)>2},  if $G$ is IBIS, then $H \in \{ {}^3\mathrm{D}_4(2):3, 2^5.\mathrm{L}_5(2)\}$ and $b(G)=3$. In each case,  straight line programs for the generators of $H$ in terms of standard generators of $G$ are available in the Online Atlas.  
We proceed by applying (T3) in the same fashion as in Section \ref{Ly}.
More explicitly, if $H={}^3\mathrm{D}_4(2):3$, then we apply (T3) by considering the centraliser of an involution $g_2$ lying in the class labelled 2A of $H$. If instead $H=2^5.\mathrm{L}_5(2)$, we again consider the centraliser of a 2A involution in $H$.  
\subsection{$\bm{G=\text{M}}$}
By Proposition \ref{b(G)>2}, it remains to check that the action of $G$ on cosets of $H=2.\text{B}$, where $b(G)=3$. 
First fix the identity coset, whose stabiliser is $H$. If $G$ is IBIS, then every irredundant base of $H$ on $G/H$ has size 2. Let $g_{17}$ be an element of order 17. We consider $K=C_H(g_{17})$ with the view of employing (T3).  We may construct $C=|C_G(g_{17})|$ in {\sf{GAP}} as a subgroup of its maximal overgroup $(L_3(2)\times S_4(2):2).2< G$ using the Online Atlas \cite[p.234]{ATLAS}. We know from the character table of $H$  \cite[p.210--219]{ATLAS} that $|K| = 136$. Now $C$ has a single conjugacy class of subgroups of order 136, so without loss of generality, we may fix one of them and set it to be $K$.  Now $C$ is divisible by 3 and 7, while $|K|$ is not, so take $g_3, g_7\in C_G(g_{17})$ to be elements of orders 3 and 7 respectively.  Clearly $H \neq Hg_3, Hg_7$, and we must also verify that $g_3g_7^{-1} \notin H$.  We observe that $k\in K$ fixes $Hg_3$ if and only if $k^{g_3} \in H$.  Since $k^{g_3}\in C$, we must have $k^{g_3} \in H\cap C =  K$, and a similar argument holds for $Hg_7$. Therefore, we check directly in our construction of $C$ and $K$ in {\sf{GAP}} that there exist $g_3$, $g_7$ such that: (1) $g_3g_7^{-1} \notin H$ so that $Hg_3$, $Hg_7$ are distinct, (2) $K_1=K^{g_3}\cap K < K$ and (3) $1< K^{g_7}\cap K_1 < K_1$. Together, these demonstrate that $K$ has a partial irredundant base of size 2 on $G/H$, hence so does $H$. Therefore, $G$ is not IBIS.

%
%

\subsection{$\bm{G=\text{B}}$}
By Proposition \ref{b(G)>2}, we need to consider
 \[H\in \{2.{}^2.E_6(2).2, 2^{1+22}.\mathrm{Co}_2, \mathrm{Fi}_{23}, 2^{9+16}.\mathrm{PSp}_8(2), \rm{Th}, (2^2\times F_4(2)):2, 2^{2+10+20}.(M_{22}:2 \times S_3)\}\]
 The action of $G$ on cosets of $H$ has base size 3 in all cases except $H=2.{}^2.E_6(2).2$, where $b(G)=4$.

%

 First let $H=2.{}^2.E_6(2).2$. We can construct $H$ as the centraliser of an element $t_0$ in conjugacy class $2A$ of $G$. The orbit lengths of $H$ on 2A involutions $t\in B$ are given in \cite[p. 216]{ATLAS} and we reproduce them in Table \ref{B_M1_orbs}.
 \begin{table}[h!]
 \begin{tabular}{ccc}
 \toprule
 $i$& Class of $t_0t_i$ & $|C_H(t_i)|$\\
 \midrule
0&1A & $2^{38}.3^9.5^2.7^2.11.13.17.19$\\ 
 1&2B & $2^{38}.3^6.5.7.11$\\
 2&2C & $2^{26}.3^6.5^2.7^2.13.17$\\
 3&3A & $2^{18}.3^9.5^2.7.11.13$\\
 4&4B & $2^{30}.3^6.5.7$\\
 \bottomrule
 \end{tabular}
 \caption{Orbit lengths of $H$ on $2A$ involutions of $B$. \label{B_M1_orbs}}
 \end{table}

Let $t_i$ be an orbit representative of the orbit labelled $i$ in Table \ref{B_M1_orbs}. We claim that $\{t_0,t_1,t_2, t_3\}$ forms a partial irredundant base for $G$, implying $G$ is not IBIS. Clearly $G_{t_0} = H$ by definition and $|G_{t_0,t_1}|  = |C_H(t_1)| = 2^{38}.3^6.5.7.11$ by Table \ref{B_M1_orbs}.  Choosing $t_2$ and $t_3$ appropriately, we find that the Sylow 2-subgroup of $|G_{t_0,t_1,t_2}|$ has order $2^{26}$, while the Sylow 2-subgroup of $|G_{t_0,t_1,t_2,t_{3}}|$ has order $2^{18}$. Therefore, these groups are distinct and non-trivial, so the claim is proved.


Now let $H=2^{1+22}.\mathrm{Co}_2$. Similar to the last case, $H$ can be constructed as the centraliser of an involution $t_0$ in class 2B of $B$, and we can instead consider the action of $G$ on its class 2B. M\"uller \cite[Table 1]{MR2379937} computed the orbit lengths of $H$ on 2B involutions.  The smallest non-trivial orbit has length 93150 and so the size of a stabiliser $K$ in $H$ of an element $t_1$ of this orbit is larger than the number of cosets of $H$ in $G$. Therefore, it is impossible to construct an irredundant base of size 3 that begins with $\{t_0,t_1\}$, so $G$ is not IBIS.

Now let $H=\rm{Fi}_{23}$. The $H$-orbits on cosets of $H$ were computed by M\"uller et al. \cite[Table 2]{MR2331753}.  The Sylow 7-subgroup of $H$ has order 7, as do the stabilisers $K_2 \cong O_8^+(3):2_2$  and $K_3\cong \mathrm{PSp}_8(2)$ of elements in the non-trivial orbits labelled 2 and 3 in \cite[Table 2]{MR2331753} respectively. Hence, taking appropriate representatives $Hx_2, Hx_3$ of these orbits, we achieve an irredundant partial base $\{H,Hx_2,Hx_3\}$ of size 3 stabilised by an element of order 7.  Since the stabiliser is not trivial, $G$ is not IBIS.

Next suppose $H = 2^{2+10+20}(M_{22}:2\times S_3)$.
Neunh\"offer et al.  \cite[Table 1]{MR2824523} computed a subset of the $H$-suborbits of $B$ acting on cosets of $H$.  Although only the orbit lengths and sizes of the corresponding stabilisers in $H$ are available, this is sufficient to show $G$ is not IBIS.  There is a single class of elements of order 5 in $H$, and two distinct orbits of $H$ on cosets with stabiliser sizes 240 and 120. Therefore,  by choosing representatives of these orbits appropriately, we can construct a partial irredundant base of size three stabilised by an element of $G$ of order 5. Therefore, $G$ is not IBIS.


For the remaining maximal subgroups $H_4=2^{9+16}.\mathrm{PSp}_8(2)$, $H_5=\rm{Th}$ and $H_6= (2^2\times F_4(2)):2$
we adopt a different approach.  In the following construction, we will use $H$ to refer to any one of the groups $H_4,H_5$ or $H_6$. Let $g_6 \in H$ have order 6, and let $g_2=g_6^3$.  We will construct distinct cosets $Hx$ and $Hy$ such that $Hx$ is fixed by $g_2$ and $g_6$, while $Hy$ is fixed by $g_2$, but not $g_6$. This will show that $\{H,Hx,Hy\}$ form a partial irredundant base for $G$ in its coset action on $H$, and so $G$ is not IBIS.

Now,  $g$ fixes $Hx$ if and only if $g^{x^{-1}}=xgx^{-1} \in H$, and similarly for $Hy$. Hence, it is sufficient to exhibit $x$ and $y$ such that $g_6^{x^{-1}}, g_2^{x^{-1}}, g_2^{y^{-1}}\in H$, while $g_6^{y^{-1}} \notin H$.  Notice that we may simplify the problem slightly by taking $x,y\in HC_G(g_2)$. We must then check that $g_6^{x^{-1}} \in H$ $g_6^{y^{-1}} \notin H$ and that $x,y, xy^{-1} \notin H$ to ensure that $H$, $Hx$ and $Hy$ are distinct.  Since $HC_G(g_6) \subseteq HC_G(g_2)$, the probability of finding $x \in C_G(g_2)$ such that $g_6^{x^{-1}} \in H$ is proportional to 
\[
\frac{|C_G(g_6)| }{ |C_G(g_2)|} = \frac{|C_G(g_6)|{|C_H(g_2)|}}{|C_G(g_2)|{|C_H(g_6)|}}.
\]
We try to maximise this quantity in our choice of the class of $g_6$, and using the Character Table Library in GAP (including the labelling of classes given there), we choose $g_6$ from classes $6C$, $6B$ and $6Q$ in $H_4$, $H_5$ and $H_6$ respectively.

We now turn to a computational search in {\sc Magma}. The relevant code to reproduce the procedure we describe is available on the author's website. We are able to construct $G$ as a matrix group in a 4370-dimensional representation over $\mathbb{F}_2$, and construct each of the maximal subgroups $H$ using straight-line programs on the standard generators, which are available in \cite{onlineATLAS}. We observe that each maximal subgroup acts reducibly on the underlying vector space $V$, so we compute the submodules preserved by each maximal subgroup.  In each case, we are able to find a relatively small submodule on which $H$ acts faithfully. We project onto this submodule so that we are able to find the conjugacy classes of $H$ using the \texttt{LMGClasses} function.  Once we have found an appropriate element $g_6$, we then find the preimage in the 4370-dimensional representation of $H$, and define $g_2=g_6^3$.  We compute $C_G(g_2)$ using \texttt{CentraliserOfInvolution} function, and search for $x$ and $y$ satisfying the above conditions by random selection.  Since determining membership of a random element of $G$ in $H$ is very time-consuming, we instead equivalently check whether the element preserves a submodule fixed by $H$. We are able to find an appropriate $x$ and $y$ in each case, so infer that $G$ is not IBIS in its coset actions on $H_4$, $H_5$ and $H_6$.


%

The proof of Theorem \ref{mainthm} is now complete.
\section*{Acknowledgements}
The author would like to thank Eamonn O'Brien for helpful discussions and assistance with the intensive computations required for the proof for the Baby Monster. The author also acknowledges the support of an Australian Research Council Discovery Early Career Researcher Award (project number DE230100579).

\bibliographystyle{abbrv}
\bibliography{ibis.bib}

\end{document}